\documentclass[hidelinks,onefignum,onetabnum]{siamart220329}

\usepackage{pgfplots} % <- Preamble

\newtheorem{Rem}{REMARK}

\newtheorem{Alg}{Algorithm}

\usepackage{lipsum}
\usepackage{amsfonts}
\usepackage{graphicx}
\usepackage{epstopdf}
\usepackage{algorithmic}
\ifpdf
  \DeclareGraphicsExtensions{.eps,.pdf,.png,.jpg}
\else
  \DeclareGraphicsExtensions{.eps}
\fi

\newsiamremark{remark}{Remark}
\newsiamremark{hypothesis}{Hypothesis}
\crefname{hypothesis}{Hypothesis}{Hypotheses}
\newsiamthm{claim}{Claim}

\headers{Trigonometric Interpolation Based Optimization}{Xiaorong Zou}
\title{Trigonometric Interpolation Based Optimization \\ for Second Order Non-Linear ODE  with Mixed Boundary Conditions
}

\author{Xiaorong Zou\thanks{Global Market Risk Analytic,  Bank of America,
  (\email{xiaorzou@gmail.com}).}}

\usepackage{amsopn}

\ifpdf
\hypersetup{
	pdftitle={On Trigonometric Approximation and Its Applications: high order ODE},
	pdfauthor={Xiaorong Zou}
}
\fi

\begin{document}
	
\maketitle
% REQUIRED
%\include{section_abstract_20241215_part2}
\begin{abstract}
In this paper, we propose a trigonometric-interpolation approach for solutions of second order nonlinear ODEs with mixed boundary conditions. The method interpolates secondary derivative $y''$ of a target solution $y$ by a trigonometric polynomial. The solution is identified through an optimization process to capture the dynamics of $y,y',y''$ characterized by the underlying differential equation.  The gradient function of the optimization can be carried out by Fast Fourier Transformation and high-degree accuracy can be achieved effectively by increasing interpolation grid points.  In case that solution of ODE system is not unique, the algorithm has flexibility to approach to a desired solution that meets certain requirements such as being positive.  Numerical tests have been conducted under various boundary conditions with expected performance. 

The algorithm can be extended for nonlinear ODE of a general order $k$ although implementation complexity will increase as $k$ gets larger. 
\end{abstract}


% REQUIRED

\section{Introduction}\label{sec:intro}
A new trigonometric interpolation algorithm was recently introduced in \cite{zou_tri_I} and \cite{zou_tri_II}.  It leverages Fast Fourier Transformation (FFT) to achieve optimal computational efficiency and converges at speed aligned with smoothness of underlying function. In addition, it can be used to approximate nonperiodic functions defined on bounded intervals.  Considering the analytic attractiveness of trigonometric polynomial,  especially in handling differential and integral operations, the proposed trigonometric estimation of a general function is expected to be used in a wide spectrum. As an example, a trigonometric interpolation based optimization method was used to solve nonlinear ordinary differential equation (ODE) of degree $1$ and the test results show that it outperforms standard Runge-Kutta scheme significantly \cite{zou_tri_II}. In this paper, we continue on applications of the trigonometric interpolation algorithm to solve the following nonlinear ODE system:

\begin{equation} \label{eq:nonlinear_ode_order2}
	y''(x) =f(x,y,y'),  \quad x\in [s,e] 
\end{equation}
\begin{eqnarray}
d_{11}y(s) + d_{12}y'(s) +d_{13}y(e) + d_{14}y'(e) &=& \alpha , \label{eq:nonlinear_ode_order2:diri}\\ % Dirichlet boundary condition 
d_{21}y(s) + d_{22}y'(s) +d_{23}y(e) + d_{24}y'(e) &=& \beta, \label{eq:nonlinear_ode_order2:neum} %Neumann boundary condition
\end{eqnarray}
where $f(x,v,u)$ is continuously differential on the range $[s, e]\times R^2$, the rank of matrix $D:=(d_{ij})_{1\le i\le 2, 1\le j\le 4}$ is $2$, and $\alpha, \beta$ are two real numbers. 

A wide variety of natural phenomena are modeled by second order ODEs. They have been applied to many problems in physics, engineering, mechanics, modern control theory, and so on. In general, closed form solutions are not available, especially in nonlinear cases, and many numerical algorithms have been proposed, often depending on the complexity of ODE (\ref{eq:nonlinear_ode_order2}) and boundary conditions (\ref{eq:nonlinear_ode_order2:diri}- \ref{eq:nonlinear_ode_order2:neum}).  

For Neumann condition with $y(s)= \alpha, y'(s)=\beta$,  it is well-know that the solution exists and is unique if  $f$ is continuous on the domain 
\begin{equation}\label{eq:dom}
Dom=\{(x, y, z)| s\le  x\le e;  (y,z)\in R^2 \}
\end{equation}
and satisfies a uniform Lipschitz
condition with respect to $y$ and $z$:
\begin{equation}\label{eq:lipschitz}
|f(x,y_1,z_1)-f(x,y_2,z_2)| \le K|y_1-y_2| + L|z_1-z_2|,
\end{equation}
and classic Runge-Kutta scheme, labeled as RK4 hereafter, provides an approximation with accuracy of order $4$ \cite{HBKeller}. The situation gets much complicated with Dirichlet condition, i.e. the boundary value problem (BVP):
\begin{equation}\label{eq:dirichlet_condi}
y(s)= \alpha, \qquad y(e)=\beta.
\end{equation}
A sufficient condition for the existence and uniqueness of a solution is established for second order linear differential equation with Dirichlet conditions in \cite{ode_naga}. As a  wide used results, we have \cite{HBKeller}
\begin{theorem}\label{thm:dirichlet_standard} 
Assume $f$ is continuous on the domain (\ref{eq:dom}) and satisfies a uniform Lipschitz condition (\ref{eq:lipschitz}).  In addition, 
$\frac{\partial f}{\partial y}$ and $\frac{\partial f}{\partial z}$ are continuous. Also,
	\begin{equation}
		\frac{\partial f}{\partial y}>0, \quad |\frac{\partial f}{\partial z}| < M,
	\end{equation}
	for some constant $M$, then BVP (\ref{eq:nonlinear_ode_order2},\ref{eq:dirichlet_condi}) has a unique solution . 
\end{theorem}
 There are quite rich researches on algorithms of numerical solutions of boundary problems of second order ODE.  Among them is the shooting method, which we apply as a benchmark model in Section \ref{sec:performance}. The method finds the proper value $\gamma$  of $y'(s)$ such that the associated solution $y(x)$ under Neumann condition $y(s)=\alpha, y'(s)=\gamma$ meets the desired boundary condition $y(e)=\beta$.  Some well-know numerical algorithms are reviewed in \cite{ode_cuma_2},  including monotone iterative methods, quasilinearization method, and finite difference method.  A Taylor approximation algorithm for solving BVP of linear ODEs was proposed in \cite{losite}. Homotopy perturbation method was introduced in \cite{ode_17_chsa} for solving BVP of linear and nonlinear  ODEs. Adomian decomposition method represents a solution in the form of Adomian polynomials and has advantage to provide analytical approximation to a wide class of nonlinear equations \cite{ode_Adomian}. 
%In \cite{ode_yuyuse},  a practical Chebyshev matrix method is proposed for the approximate solution of high order linear Fredholm Integro-differential equations (FIDES) with constant coefficients under the initial-boundary conditions in terms of Chebyshev polynomials. This method transforms the integro-differential equation to a system of linear algebraic equations. It especially fits for direct solution of constant coefficients FIDEs. 
In \cite{ode_19_gecu}, an algorithm based on a combination of the Adomian decomposition method (ADM) and the reproducing kernel method (RKM) is proposed for solving nonlinear BVPs. Cubic spline method is used to solve a class of singular two-point boundary value problems of second order linear ODE in \cite{ode_24_kare}. The reproducing kernel method has developed rapidly in solving BVPs of high order linear ODE in recent years. In \cite{ode_zhlish}, a new multiscale algorithm was proposed to solve BVPs of second order linear ODE. Some other references in that direction can be found in \cite{ode_8}, \cite{ode_bumozh13}, \cite{ode_3}, \cite{ode_22}, and  \cite{ode_36}. In general,  certain reproducing kernel is used to construct orthogonal basis for a target Hilbert space and the solution of ODE problem can be estimated by the projection to a finite subspace.  The approximation can be obtained by solving linear system based on the boundness of the associated linear differential differential operator of ODE. 

Mixed Dirichlet and Neumann boundary (\ref{eq:nonlinear_ode_order2:diri}- \ref{eq:nonlinear_ode_order2:neum}) has been considered in \cite{losite} with linear second ODE. The solution is represented by a multi-point Taylor expansion. %The number and location of the base points of that expansion are carefully chosen to guarantee the convergence.

In this paper, we aim to extend the trigonometric interpolation based optimization algorithm, developed in \cite{zou_tri_II}, for the solutions of second order ODE  (\ref{eq:nonlinear_ode_order2}-\ref{eq:nonlinear_ode_order2:neum}). The proposed algorithm, labeled as TIBO hereafter, represents the secondary derivative $y''$ of a target solution $y$ by a trigonometric interpolation polynomial. The boundary conditions can be captured by two parameters in the derived close-form expression of $y$, similar to how initial condition is captured in Adomian decomposition method. More importantly, the values of derivatives $\{y'(x_k), y''(x_k)\}_{k\le M}$  can be effectively carried out through $\{y(x_k)\}_{k\le M}$ at grid points $\{x_k\}_{k\le M}$ of an interpolation. Further more,  FFT can be leveraged to optimize the required numerical computation such that large number of grid points can be applied to improve the accuracy of the approximation.  The combination of these features make TIBO an effective global optimization algorithm for general ODE systems.  In fact, it can be extended for nonlinear ODE of a general order $K$ although implementation complexity will increase as $K$ gets larger.  

TIBO is expected to have high accuracy when underlying target function $y$ is sufficient smooth as stated in Theorem \ref{conv_k_new}, which is confirmed by the numerical tests conducted in Section \ref{sec:performance}.  Another advantage of global optimization is its flexibility to integrate some extra conditions into optimization so it converges to a desired solution in case of multiple solutions.  In Section \ref{subsec:Dirichlet},  we provide an example with two solutions and TIBO with conditioning on a range of $y'(s)$ converges to the right solution.  In Section \ref{subsec:Mix_1},  condition on lower bound of solution is used to identify the desired solution from ODE system with two solutions.   Note finding  positive solutions of an ODE system is an interesting topic \cite{ode_agwoli}  and \cite{ode_wo1}.

TIBO is comparable to Adomian decomposition method  in the sense that it uses analytical close form to globally approximate the solutions and avoid errors due to discretization or Taylor-formula-based local approximation. Using trigonometrical polynomial and applying optimization process, TIBO is flexible to handle nonlinear terms and leverages FFT to apply large number of interpolation points to increase accuracy, thus is expected to outperforms polynomial-based Adomian decomposition method.

The rest of paper is organized as follows.  In Section \ref{sec:trig}, we summarize the relevant results of trigonometric interpolation algorithm developed in \cite{zou_tri_I} and \cite{zou_tri_II}.  Section \ref{sec:algorithm} is devoted to develop TIBO. The main idea is similar to what is used to address first order ODE in \cite{zou_tri_II} with extra attention to copy with more complicated objective function and initial conditions.  A compressive numerical tests are conducted in Section \ref{sec:performance} for performance assessment with three types of boundary conditions.  We compare TIBO with a benchmark method that combines shooting method and classic Runge-Kutta scheme \cite{ryts} and report the results in three subsections \ref{subsec:Neumann}-\ref{subsec:Mix_1}.  The summary is made on Section \ref{sec:summary}.

\section{Trigonometrical Interpolation on Non-Periodic Functions}\label{sec:trig}
In this section \ref{sec:trig}, we review the relevant results of trigonometric interpolation algorithm developed in \cite{zou_tri_I} and \cite{zou_tri_II} starting with the following interpolation algorithm. 
\begin{Alg}\label{main_thm} Let $f(x)$ be an odd periodic function \footnote{Similar results for even periodic function is also available in \cite{zou_tri_I}.} with period $2b$ and $N=2M=2^{q+1}$ for some integer $q\ge 1$ and $x_j, y_j$ are defined by	
	\begin{eqnarray}
		x_j &:=& -b + j\lambda,  \quad \lambda= \frac{2b}N , \quad 0\le j <N, \label{x_grid_N} \\
		y_j &:=& f(x_j),  \label{f_N_interpolation_new} 
	\end{eqnarray}
	then there is a unique $M-1$ degree trigonometric polynomial
		\begin{eqnarray*}
			f_M(x) &=& \sum_{0\le j <M}a_j \sin\frac{j\pi x}b, \label{f_M_interpolation_new_odd}\\
			a_j&=&\frac{2}N\sum_{0\le k <N} (-1)^j y_k \sin\frac{2\pi j k}{N}, \quad 0\le j <M \label{aj_odd}
		\end{eqnarray*}
		such that it fits to all grid points, i.e.
		\begin{equation*}\label{error_even_new_odd} 
			f_M(x_{k})=y_{k}, \quad  0\le k <N.
		\end{equation*}
\end{Alg} 
One can computer coefficients by Inverse Fast Fourier Transform (ifft):
\[	
\{a_j (-1)^j\}_0^{N-1} = 2\times Imag (ifft(\{y_k\}_{k=0}^{N-1}))
\]
%Note that the uniqueness of solution provides the theoretic support for the algorithm proposed in Section \ref{sec:algorithm}.
It is shown in \cite{zou_tri_I} that  $f_M(x)$ and associated derivatives converge in desired order. 
\begin{theorem}\label{conv_k_new}
Let $f(x)$ be an periodic function with period $2b$ and
 $|f^{(K+1)}(x)|$ exist with an upper bound $D_{K+1}$, then 
	\begin{eqnarray}
		|f_M(x) - f(x)| &\le& \frac{C_1({D_{K+1}})}{N^{K}},   \label{detal_M_f} \\
		|f^{(k)}_M(x) - f^{(k)}(x)| &\le& \frac{C_2({D_{K+1}})}{N^{K-k}}, \quad 1\le k <K. \label{der_detal_M_f}
	\end{eqnarray}
	where $C_1({D_{K+1}})$ and $C_2({D_{K+1}})$ are two constants depending on $D_{K+1}$.
\end{theorem}	
In \cite{zou_tri_II}, the algorithm \ref{main_thm} has been enhanced so it can be applied to a nonperiodic function $f$ over a bounded interval $[s,e]$ whose $K+1$-th derivative $f^{(K+1)}(x)$ exists. To seek for a periodic extension with same smoothness, we assume that $f$ can be extended smoothly such that $f^{(K+1)}$ exists and is bounded over $[s-\delta, e+\delta]$ for certain $\delta>0$.  A smooth periodic extension of $f$ can be achieved by a cut-off smooth function $h(x)$ with following property: 
%One can use the following method to periodically extend $g$ by a cut-off function $h$ with the following properties 
\begin{eqnarray*}
	h(x)=\left\{\begin{array}{cc}
		1 & x\in [s,e], \\
		0 & \mbox{$x<s-\delta$ or $x>e+\delta$}. \\
	\end{array}\right.
\end{eqnarray*}
A cut-off function with closed-form analytic expression is proposed  in \cite{zou_tri_II}.  Let
	\begin{equation}\label{ob}
	o=s-\delta, \quad b=e+\delta -o,
\end{equation}
and define $F(x):=h(x+o)f(x+o)$ for $x\in [0,b]$. One can treat $F(x)$ as an odd periodic function with period $2b$. Apply Algorithm \ref{main_thm} to generate the trigonometric interpolation of degree $M$ with $N=2M$ evenly-spaced grid points over $[-b,b]$
\[
	F_M(x) = \sum_{0\le j < M} a_j \sin\frac{j\pi x}{b},
\]
and let 
\begin{equation*}\label{fMext}
	\hat{f}_{M}(x)=F_M(x-o)=\sum_{0\le j < M} a_j \sin\frac{j\pi (x-o)}{b}. 
\end{equation*}
 $\hat{f}_{M}(x)|_{[s,e]}$ can be treated as an trigonometric interpolation of $f$ since $\hat{f}_{M}(x_k)=f(x_k)$ for all grid points $x_k\in [s,e]$.  Numerical tests on certain basic functions demonstrates that $\hat{f}$ approach $f$  with decent accuracy when $f$ is sufficient smooth \cite{zou_tri_II}. 
 
\section{The development of TIBO}\label{sec:algorithm}
In this section, we aim to develop TIBO, the algorithm mentioned in Section \ref{sec:intro}, to solve ODE system (\ref{eq:nonlinear_ode_order2}-\ref{eq:nonlinear_ode_order2:neum}).  We assume that $f(x,v,u)$ in Eq (\ref{eq:nonlinear_ode_order2}) is continuous differential on $[s-\delta, e+\delta]\times R^2$ for certain $\delta>0$.  By parallel shifting if needed,  we assume $s-\delta=0$ without loss of generality. Let $h$  be a cut-off function specified in Section \ref{sec:trig} and extend $f(x,v,u)$ to $F(x,v,u)$ as follows
\[
	F(x,v, u)=f(x,v, u)h(x)  \qquad x \in [0,b]
\]
Consider a solution $v(x)$ of the following ODE system
\begin{eqnarray}
	v''(x) &=& F(x,v(x),v'(x)),  \quad x\in [0,b] \label{eq:nonlinear_ode_order2_F},\\
	\alpha  &=& d_{11}v(s) + d_{12}v'(s) +d_{13}v(e) + d_{14}v'(e), \label{eq:nonlinear_ode_order2:diri_F}\\ % Dirichlet boundary condition 
	\beta  &=& d_{21}v(s) + d_{22}v'(s) +d_{23}v(e) + d_{24}v'(e). \label{eq:nonlinear_ode_order2:neum_F} 
	%v(s) &=&\alpha, \qquad v(e)=\beta , \label{eq:nonlinear_ode_order2:diri_F}\\ % Dirichlet boundary condition 
	%u(s) &=&\alpha, \qquad v(s)=\beta , \label{eq:nonlinear_ode_order2:neum_F} 
\end{eqnarray}
It is clear that $v(x)|_{[s,e]}$ solves ODE (\ref{eq:nonlinear_ode_order2}-\ref{eq:nonlinear_ode_order2:neum}). Let $u(x)=v'(x)$ and $z(x)=v''(x)$ to simplify the notations. 
It is clear that $z(x)$ vanishes at boundary point $0,b$ as well as its derivatives $z^{(k)}$, hence it can be smoothly extended as an odd periodic function with period $2b$, which can be approximated by odd trigonometric polynomial using Algorithm \ref{main_thm}.  To that direction, let  $\{(x_k, z_k)\}_{k=0}^{N-1}$ be a grid set of $z(x)$ such that
	\begin{eqnarray*}
		x_k &=& -b + \frac{2b }N k, \qquad k=0, 1, \cdots, N-1, \label{eq:grid_x_ode}\\
		z_0&=&0, \quad z_k = -z_{N-k},  \qquad 1\le k <M. \label{eq:grid_z_ode}
	\end{eqnarray*}
Note 
\[
s = x_{M+m}, \qquad e=x_{M+m+n}.
\]
Assume 
\begin{equation}\label{eq:z_M}
	z_M(x) = \sum_{0\le j<M}b_j \sin \frac{j\pi x}{b} 
\end{equation}
is the interpolant of $z(x)$ by Algorithm \ref{main_thm} such that 
\begin{equation}\label{eq:b_j_ext}
b_j=\frac{2}N\sum_{k=0}^{N-1} (-1)^j z_k \sin\frac{2\pi j k}{N}=\frac{4}N\sum_{k=M}^{N-1} (-1)^j z_k \sin\frac{2\pi j k}{N}, \quad 0\le j < M. 
\end{equation}
$u$ and $v$ can be derived accordingly 
\begin{eqnarray}
	\tilde{u}_M(x ) &=& a_0 -   \frac{b}{\pi }\sum_{1\le j <M} \frac{b_j}{j} \cos  \frac{j\pi x}{b},  \label{eq:tildeu}\\
	\tilde{v}_M(x ) &=& a_1 +  a_0 x -  (\frac{b}{\pi })^2 \sum_{1\le j <M} \frac{b_j}{j^2}  \sin  \frac{j\pi x}{b},  \label{eq:tildev}
\end{eqnarray}
We shall adapt the following notations in the rest of paper for convenience. 
\begin{eqnarray*}
	u_k &=& \tilde{u}_M(x_k), \quad v_k = \tilde{v}_M(x_k), \\
	F_k &=& F(x_k, v_k, u_k), \quad DF_{k;u} =\frac {\partial F}{\partial u} (x_{k}, v_k, u_{k}),\quad DF_{k;v} =\frac {\partial F}{\partial v} (x_{k}, v_k, u_{k})\\
	X &=& \{x_k\}_{0\le k<N}, \quad  Z= \{z_k\}_{0\le k<N}, \quad U= \{u_k\}_{0\le k<N}, \quad V= \{v_k\}_{0\le k<N}\\
	F &=& \{F_k\}_{0\le k<N}, \quad  DF_u = \{DF_{k;u}\}_{0\le k<N},\quad DF_v = \{DF_{k;v}\}_{0\le k<N},\\
	X^R &=& \{x_k\}_{M\le k<N}, \quad  Z^R= \{z_k\}_{M\le k<N}, \quad U^R= \{u_k\}_{M\le k<N}, \quad V^R= \{v_k\}_{M\le k<N}\\
	F^R &=& \{F_k\}_{M\le k<N}, \quad  DF^R_u = \{DF_{k;u}\}_{M\le k<N},\quad DF^R_v = \{DF_{k;v}\}_{M\le k<N}.
\end{eqnarray*}
ODE system (\ref{eq:nonlinear_ode_order2_F}-\ref{eq:nonlinear_ode_order2:neum_F}) can be solved by minimizing the following objective function:  
\begin{equation}\label{target_phi}
	\phi(z_{M}, z_{M+1}, ...,z_{N-1}) = \frac1{2M}\sum_{M\le k < N} (z_k-F_k)^2.  
\end{equation}
\begin{Rem}
	  $z_M(x)$ is uniquely determined by $\{z(x_k)\}$ by Algorithm \ref{main_thm} for a given $z$.  As such,  it is expected that  $\phi(z_{M}, z_{M+1}, ...,z_{N-1})$ converges to $0$ if $z$ does exists. 
\end{Rem}
We need an effective way to calculate its gradient $\frac{\partial \phi}{\partial Z}$ when $M$, the number of variables of $\phi$, is not small. For $M\le t <N$,
\begin{equation}\label{gradient_v2}
	M\frac {\partial \phi}{\partial z_t} = (z_t-F_t) - \sum_{M\le k<N} (z_k-F_k) DF_{k,u} \frac {\partial u_k}{\partial z_t} - \sum_{M\le k<N} (z_k-F_k) DF_{k,v} \frac {\partial v_k}{\partial z_t}.
\end{equation}
It is clear 
\begin{equation}\label{b_j_z_k}
	\frac{\partial b_j}{\partial z_t} = \frac{4}N (-1)^j  \sin\frac{2\pi j t}{N}, \quad 0\le j< M, \quad M\le t< N. 
\end{equation}
Let $0_M$ be the $M$-dim zero vector and define $\frac{1}0:=0$,  $sum(T)$ be the summation of all elements in a given vector $T$ and 
$A_1\circ A_2$ denote the Hadamard product that applies the element-wise product to two metrics of same dimension. Adapt the following notations to copy with boundary conditions.
\begin{eqnarray*}
	J &=& [0, 1, \frac12, \dots, \frac1{M-1}, 0_M],  \qquad A = [(-1)^{j}]_{0\le j <N}\\
	%\Phi&=&[\sin \frac{2\pi j (m+n) }N ]_{j=0}^{M-1}, \quad \Phi_{N} = [\Phi, 0_M], \\
	\Phi^{s}_m&=&[ \sin \frac{2\pi j m }N ]_{j=0}^{M-1}, \quad \Phi^s_{m,N} = [\Phi^s, 0_M], \\
	\Phi^{s}_{(m+n)}&=&[ \sin \frac{2\pi j {m+n} }N ]_{j=0}^{M-1}, \quad \Phi^s_{m+n,N} = [\Phi^s_m, 0_M], \\
	\Phi^{c}_m&=&[ \cos \frac{2\pi j m }N ]_{j=0}^{M-1}, \quad \Phi^c_{m, N} = [\Phi^c_m, 0_M], \\
	\Phi^{c}_{(m+n)}&=&[ \cos \frac{2\pi j {m+n} }N ]_{j=0}^{M-1}, \quad \Phi^c_{m+n,N} = [\Phi^c_{m+n}, 0_M], \\
	\Psi_u&=&\{ (z_j-F_j)DF_{j;u} \}_{j=M}^{N-1},  \quad \Psi^u_N = [0_M, \Psi_u],  \quad  I_u = sum(\Psi_u)\\
	\Psi_v&=&\{ (z_j-F_j)DF_{j;v} \}_{j=M}^{N-1},  \quad \Psi^v_N = [0_M, \Psi_v],  \quad  I_v = sum(\Psi_v), II_v = sum(\Psi_v\circ X^R),\\
S_m &=& \frac{b^2}{\pi^2}\sum_{0<j<M} \frac{b_j}{j^2} \sin\frac{2\pi jm}{N}, \quad S_{m+n} =\frac{b^2}{\pi^2} \sum_{0<j<M} \frac{b_j}{j^2} \sin\frac{2\pi j(m+n)}{N} \\
C_m &=& \frac{b}{\pi }  \sum_{0<j<M} \frac{b_j}{j} \cos\frac{2\pi jm}{N} , \quad C_{m+n} = \frac{b}{\pi }  \sum_{0<j<M} \frac{b_j}{j} \cos\frac{2\pi j(m+n)}{N}.
\end{eqnarray*}
By Eq. (\ref{eq:tildeu}) and (\ref{eq:tildev}), we have
\begin{eqnarray*}
	v(s) &=& v(x_{M+m}) = a_1 + a_0s - S_m, \label{eq:init_vs},\\
	v(e) &=& v(x_{M+m+n}) = a_1 + a_0e - S_{m+n}, \label{eq:init_ve} \\
	u(s) &=& u(x_{M+m}) =  a_0 -C_m, \label{eq:init_us} \\
	u(e) &=& u(x_{M+m+n})=a_0 - C_{m+n}. \label{eq:init_ue}
\end{eqnarray*}
Let $\nabla w:=(\frac {\partial w}{\partial z_i})_{M\le i <N}$ be the gradient vector for a variable $w=w(z_M, \dots z_{N-1})$. Applying Eq (\ref{b_j_z_k}), we obtain
\begin{eqnarray}
	\nabla S_m &=& \frac{4b^2}{\pi^2} Imag (ifft(A \circ \Phi^s_{m,N}\circ J^2)) \label{eq:nabla_s_m},\\
	\nabla S_{m+n}  &=& \frac{4b^2}{\pi^2} Imag(ifft(A \circ \Phi^s_{m+n,N}\circ J)) \label{eq:nabla_s_m_n},\\
	\nabla C_m  &=& \frac{4b^2}{\pi^2} Imag(ifft(A \circ \Phi^c_{m,N}\circ J^2)) \label{eq:nabla_c_m},\\
	\nabla C_{m+n}  &=& \frac{4b^2}{\pi^2} Imag(ifft(A \circ \Phi^c_{m+n,N}\circ J)) \label{eq:nabla_c_m_n}.
\end{eqnarray} 
One can solve $a_0$ and $a_1$ by Eq (\ref{eq:nonlinear_ode_order2:diri_F}) and (\ref{eq:nonlinear_ode_order2:neum_F}), 
\begin{equation}\label{eq:a0a1}
	\begin{pmatrix}
		a_{0} \\
		a_{1} 
	\end{pmatrix} =\frac{1}{m_{11}m_{22}-m_{12}m_{21}} 
\begin{pmatrix}
	m_{22} & -m_{12} \\
	-m_{21} & m_{11} 
\end{pmatrix}
\begin{pmatrix}
	\alpha + \mu \\
	\beta  + \nu 
\end{pmatrix},
\end{equation}
where 
\begin{eqnarray*}
m_{11} &=& d_{11}s + d_{12} + d_{13}e +d_{14}, \quad m_{12} = d_{11} + d_{13},\\
m_{21} &=& d_{21}s + d_{22} + d_{23}e +d_{24}, \quad m_{22} = d_{21} + d_{23}.\\
\end{eqnarray*}
and 
\begin{eqnarray*}
	\mu &=& d_{11}S_m + d_{12}C_m + d_{13}S_{m+n} +d_{14}C_{m+n}\\
	\nu &=& d_{21}S_m + d_{22}C_m + d_{23}S_{m+n} +d_{24}C_{m+n}.
\end{eqnarray*}
We obtain by Eq (\ref{eq:a0a1}), (\ref{eq:nabla_s_m}-\ref{eq:nabla_c_m_n})
\begin{eqnarray*}
	\nabla a_0 &=&  \frac{1}{m_{11}m_{22}-m_{12}m_{21}} (m_{22} \nabla \mu  -m_{12} \nabla \nu )   \\
	\nabla a_1  &=& \frac{1}{m_{11}m_{22}-m_{12}m_{21}} (-m_{21} \nabla \mu   + m_{11} \nabla \nu )
\end{eqnarray*}
where $\nabla I$ and $\nabla II$ can be calculated by (\ref{eq:nabla_s_m}-\ref{eq:nabla_c_m_n}).
It is worthwhile to point out that $\nabla a_0$ and $\nabla a_1$ are constant vectors, which makes the optimization process effective.
We like to solve $U$ in term of $Z$ to copy with $\frac{\partial U}{\partial Z}$ in Eq. (\ref{gradient_v2}).  By Eq (\ref{eq:b_j_ext}) and $z_0=z_M=0$, we obtain for $M\le k <N$
\begin{eqnarray}
	u_k &=&  a_0 - \sum_{0\le j <M} (-1)^j \frac{bb_j}{j\pi }  \cos  \frac{2 \pi jk }{N} \nonumber\\
	&=&  a_0 - \frac{2b}{\pi N} \sum_{0\le j <M} \frac{1}{j}  \cos \frac{2 \pi jk }{N} \sum_{0\le l<N} z_l \sin\frac{2\pi j l}{N} \label{uk_1}\\
	&=& a_0 - \frac{4b }{\pi N}\sum_{M\le l<N} z_l \sum_{0\le j <M}  \frac1j \cos  \frac{2 \pi jk }{N}   \sin\frac{2\pi j l}{N}. \label{uk_2} 
\end{eqnarray}
The last step is due to $z_{l}\sin\frac{2\pi j l}{N}=z_{N-1}\sin\frac{2\pi j (N-l)}{N}$. 
One can rewrite Eq. (\ref{uk_1}) in term of ifft as follows:
\begin{equation*}\label{eq:U}
	U = a_0 - \frac{2bN}{\pi} Re\{ifft(J\circ Im\{ifft(Z)\})\}[M:N-1].
\end{equation*}
Similarly, we need express $V$ in term of $Z$.  
\begin{eqnarray}
	v_k &=&  a_1 + a_0 x_{k} - \frac{b^2}{\pi^2 }  \sum_{0\le j <M} (-1)^j \frac{b_j}{j^2 }  \sin  \frac{2 \pi jk }{N} \nonumber\\
	&=&  a_1 +a_0 x_{k}  - \frac{2b^2}{\pi^2 N} \sum_{0\le j <M} \frac{1}{j^2}  \sin \frac{2 \pi jk }{N} \sum_{0\le l<N} z_l \sin\frac{2\pi j l}{N} \label{vk_1}\\
	&=& a_1 +a_0 x_{k} - \frac{4b^2}{\pi^2 N}\sum_{M\le l<N} z_l \sum_{0\le j <M}  \frac{1}{j^2}  \sin\frac{2 \pi jk }{N}   \sin\frac{2\pi j l}{N}. \label{vk_2} 
\end{eqnarray}
By (\ref{vk_1}),
\begin{equation}\label{eq:V}
	V^R = a_1 + a_0 X^R - \frac{2b^2N}{\pi^2} Im\{ifft(J^2\circ Im\{ifft(Z)\})\}[M:N-1],
\end{equation}
By (\ref{uk_2}) and (\ref{vk_2}), we obtain 
\begin{eqnarray}	
	\frac{\partial u_k}{\partial z_t} &=& \frac{\partial a_0}{\partial z_t} -\frac{4b }{\pi N} \sum_{0\le j <M}  \frac 1j \cos  \frac{2 \pi jk }{N}   \sin\frac{2\pi j t}{N}  \label{uz_2}, \quad M\le t<N. \\
	\frac{\partial v_k}{\partial z_t} &=& \frac{\partial a_1}{\partial z_t} + x_k\frac{\partial a_0}{\partial z_t}  -\frac{4b^2 }{\pi^2 N} \sum_{0\le j <M}  \frac 1{j^2} \sin  \frac{2 \pi jk }{N}   \sin\frac{2\pi j t}{N} \label{vz_2}, \quad M\le t<N.
\end{eqnarray}
We are ready to attack the nontrivial term in Eq. (\ref{gradient_v2}). Define 
\begin{eqnarray}
\phi^u_t &:=&\sum_{M\le k<N} (z_k-F_{k,u}) DF_{k,u} \frac {\partial u_k}{\partial z_t}, \label{eq:phi_u_t}\\
\phi^v_t &:=&\sum_{M\le k<N} (z_k-F_{k,v}) DF_{k,v} \frac {\partial v_k}{\partial z_t}. \label{eq:phi_v_t}
\end{eqnarray}
By (\ref{uz_2}) and (\ref{vz_2}),
\begin{eqnarray*}
	\phi^u_t &=&  \sum_{M\le k<N} (z_k-F_{k;u}) DF_{k;u}  \frac{\partial a_0}{\partial z_t} \\
	&-& \frac{4b }{ \pi N}  \sum_{0\le j <M} J_j \sin\frac{2\pi j t}{N}   \sum_{M\le k<N} (z_k-F_{k;u}) DF_{k;u}  \cos  \frac{2 \pi jk }{N}. \label{wuz}\\
\phi^v_t &=&  \sum_{M\le k<N} (z_k-F_{k;v}) DF_{k;v}  (\frac{\partial a_1}{\partial z_t} + x_k \frac{\partial a_0}{\partial z_t})\\
&-& \frac{4 b^2}{ \pi^2 N}  \sum_{0\le j <M} J^2_j \sin\frac{2\pi j t}{N}   \sum_{M\le k<N} (z_k-F_{k;v}) DF_{k;v}  \sin  \frac{2 \pi jk }{N}. \label{wvz}
\end{eqnarray*}
The gradient vector (\ref{gradient_v2}) can be formulated by FFT as follows:
	\begin{eqnarray}
		\phi^u &=& I_u \nabla a_0 -  \frac{4b N}{\pi} Im \{ifft (J\circ Re[ifft(\Psi^u_N)])\}[M:N-1] \label{eq:phi_u}\\
		\phi^v &=& I_v \nabla a_1  + II_v \nabla a_0 -  \frac{4b^2 N}{\pi^2} Im \{ifft (J^2 \circ Im[ifft(\Psi^v_N)])\}[M:N-1]\label{eq:phi_v}
	\end{eqnarray}
	
The gradient (\ref{gradient_v2}) is reduced to
\begin{equation}\label{eq:graident_order2}
	%\frac{\partial \phi}{\partial Z}= \frac{1}{M}(Z^R-F^R-\phi^u -\phi^v).
	\nabla \phi= \frac{1}{M}(Z^R-F^R-\phi^u -\phi^v).  
\end{equation}
In summary, we obtain a trigonometric interpolation on the solution of ODE (\ref{eq:nonlinear_ode_order2}-\ref{eq:nonlinear_ode_order2:neum}) \footnote{We assume that $s=\delta$ in Algorithm \ref{alg_tibo}. In general,  one can first parallel shift $x$ domain to left by $e-s$.}.
\begin{Alg}\label{alg_tibo}
	\begin{enumerate}
\item Construct cut-off function $h$ with parameter $s,e,\delta$.
\item Solve optimization problem (\ref{target_phi}) with gradient vector defined by Eq. (\ref{eq:graident_order2}), where $\phi^u$ and $\phi^v$ are defined by Eq. (\ref{eq:phi_u}) and (\ref{eq:phi_v}).
\item Use the optimization output vector $Z$ to find coefficients $B$ by Eq  (\ref{eq:b_j_ext}).
\item Apply $B$ to find $a_0$ and $a_1$ by Eq. (\ref{eq:a0a1}) and thus obtain the trigonometric approximation $\tilde v_M$ (\ref{eq:tildev}).		
	\end{enumerate}
\end{Alg}

\section{The numerical performance assessments}\label{sec:performance}
\subsection{The description of performance tests}\label{subsec:description}
Let $f(x; \theta) =x\cos\theta x$, $s=1, e=3$, for fixed $A=(a_{ij})_{2\times 4}$ and $C=(c_{uu},c_{uv}, c_{vv}, c_u, c_v)$, 
one can easily see that $f$ solves the following
\begin{eqnarray}
y'' &=& f''(x) - c_{uu} (f'(x))^2 - c_{uv} f(x)f'(x) - c_{vv}f^2(x) -c_u f'(x) - c_v f(x) \label{eq:example_ode} \\
&+& c_{uu} (y')^2 + c_{uv} yy' + c_{vv}y^2 + c_u y' + c_v y \nonumber\\
\alpha  &=& d_{11}y(s) + d_{12}y'(s) + d_{13}y(e) + d_{14}y'(e) \label{eq:example_ode_cond1}\\
\beta  &=& d_{21}y(s) + d_{22}y'(s) + d_{23}y(e) + d_{24}y'(e) \label{eq:example_ode_cond2}
\end{eqnarray}
where $\alpha$ and $\beta$ are defined by
\begin{eqnarray*}
\alpha  &=& d_{11}f(s) + d_{12}f'(s) + d_{13}f(e) + d_{14}f'(e)\\
\beta  &=& d_{21}f(s) + d_{22}f'(s) + d_{23}f(e) + d_{24}f'(e).
\end{eqnarray*}
 $f(x; \theta) =x\cos\theta x$ will be refereed as base solution, denoted by $y_b$ in the rest of the section. 

We shall use the following setting as default in all performance test conducted in this section unless otherwise specified. 
\begin{equation}\label{eq:C_default}
C=(0.1, 0.1, 1.0, 0.1, 1.0).
\end{equation}
The performance tests will include following boundary conditions:
\begin{table}[htbp]
	\centering
	\small
	\caption{The four testing boundary conditions }
	\begin{tabular}{lrrrrrrrrl}
	 Case	& \multicolumn{1}{l}{$d_{11}$} & \multicolumn{1}{l}{$d_{12}$} & \multicolumn{1}{l}{$d_{13}$} & \multicolumn{1}{l}{$d_{14}$} & \multicolumn{1}{l}{$d_{21}$} & \multicolumn{1}{l}{$d_{22}$} & \multicolumn{1}{l}{$d_{23}$} & \multicolumn{1}{l}{$d_{24}$} & condition on \\ \hline\hline
	 $Neumann$ & 1     & 0     & 0     & 0     & 0     & 1     & 0     & 0     & $v_s,u_s$ \\
		$Dirichlet$ & 1     & 0     & 0     & 0     & 0     & 0     & 1     & 0     & $v_s, v_e$ \\
		$Mix$   & 1     & 1     & 0     & 0     & 0     & 0     & 1     & 1     & $v_s+u_s, v_e+u_e$ \\ \hline
		%$Mix_2$ & 1     & 0     & 1     & 0     & 0     & 1     & 0     & 1     & $v_s+v_e, u_s+u_e$ \\ \hline
	\end{tabular}%
	\label{tab:test}%
\end{table}%
Note that the first two types are the classic Neumann and Dirichlet conditions. %while $Mix_2$ conditions on give values of $y(s)+y(e)$ and $y'(s)+y'(e)$.  
We use classic Runge-Kutta scheme \cite{ryts} (Section 9.4.1, page 284), labeled $rk4$, as benchmark for performance assessment on Neumann initial condition $\alpha=y(s), \beta=y'(s)$.  For Dirichlet boundary condition, we follow the shooting method, more specifically,  we search for $y'(s)$ such that $y(e)=\beta$ and then apply Runge-Kutta scheme.  For $mix$, we search pair $(y(s),y'(s))$ such that Equations  (\ref{eq:example_ode_cond1}-\ref{eq:example_ode_cond2}) are satisfied, and then apply Runge-Kutta scheme with identified $(y(s),y'(s))$.

The performance of shooting method is sensitive to the initial guess on $y'(s)$  under Dirichlet boundary condition and $y(s), y'(s)$ under $Mix$ boundary condition.   In general, TIBO's performance is expected to be sensitive to initial value of related variables, which is provided by benchmark model in our testing.  Especially, initial value can determine how optimization process converges when solution is not unique.  For each of three test types, we construct $25$ test scenarios associated to initial values created as follows: 
\begin{enumerate}
	\item Randomly simulate $10$ numbers in the range $[-0.5,0.5]$, place them in two groups as follows.   
	\begin{eqnarray*}
		init_{y} &=& [0.41,0.41,-0.40,0.05, 0.47] \\
		init_{y'} &=& [0.31,-0.37, 0.13,-0.22,0.46]
	\end{eqnarray*}
	\item Create $25$ initial $y(s), y'(s)$ by $5$ groups for $Mix$
	\[
	[f(s) + i\cdot init_{y}, f'(s) + i\cdot init_{y'}], \quad i=1, 2, -2, 3,-3.
	\]
  $init_{y}$ is adjusted to $f(s)$ under $Dirichlet$.
	\item For a pair $(i, j)$ ($1\le i, j \le 5$), the scenario in $j$ position of group $i$ is assigned with $ID=(i-1)\cdot 5 + j$. For example,  $(y(s) + 3\cdot init_{y}(1), y'(s) + 3\cdot init_{y'}(1))$ is assigned scenario id $(4-1)\cdot 5+1=16$. 
\end{enumerate}
For each of three test types, we conduct two sets of tests with $\theta\in (\frac{\pi}2, \frac{3\pi}2,)$ and each set contains $25$ scenarios described above. $f(x;\theta)$ becomes more volatile as $\theta$ gets larger and, model performance with $\theta=\pi/2$ is expected to outperform that with $\theta=3\pi/2$. For a given test type and $\theta$, overall performance on $25$ scenarios provide us how robust the optimization process is. 

We report the followings for each of three test types in Table \ref{tab:test}. 
\begin{enumerate}
	\item $max|y''-f|$:  The max difference is calculated by applying identified approximation $y_{opt}$, i.e. $y_{opt}''-f(x,y_{opt}, y'_{opt})$.  The max is taken over the set of $2^{10}$ equally-spaced points over $[0,b]$.  Note that grid point set used in the optimization algorithm is determined by $q=7$. As such, negligible $max$ indicates that $y$ converges to the solution over $[0,b]$.
	\item  $max(|y_{opt}-y_{b}|)$.  The max is taken over points over $[s,e]$ generated by $q=10$ as explained above.  A negligible value indicates that the optimization converges to the base $y_b$.  
	\item   $max(|y_{opt}-y_{s}|)$: similar max difference between $y_{opt}$ and another solution assuming that there are more than one solutions.  
	\item $max(|y_{rk4}-y_{b}|)$: the max difference between benchmark $rk4$ scheme and base solution over the points in $[s,e]$ generated by $q=10$  as explained above.
	\item $status$:  The convergent status of identified solution $y_{opt}$ by optimization.  There are three possible cases:  
	\begin{itemize}
		\item $\to y_b$: converges to base solution $y_b$ with negligible $max(|y''-f|)$ and $max(|y_{opt}-y_{b}|)$,
		\item  $\to y_s$: converges to a different solution $y_s$ with negligible $max(|y''-f|)$ but significant $max(|y_{opt}-y_{b}|)$.
		\item Does not converge to any solutions with significant  $max(|y''-f|)$ \footnote{The optimization process can either reach to a local minimal value or stop without reaching any local minimum. We treat those cases as divergent in this paper.}. 
	\end{itemize}
	\item $number$:  The number of scenarios with the specified status.
\end{enumerate}
\subsection{The performance with Neumann condition}\label{subsec:Neumann}
Table \ref{table:Neumann} summarizes the test results.  As expected,  the optimization process converges to the unique base solution $y_b$.  Note that $y_{opt}$ outperforms $y_{rk4}$ significantly in both case and performance with $\theta=\pi/2$ outperform with $\theta=\frac{3\pi}2$ as expected.
\begin{table}[htbp]
	\centering
	\small
	\caption{The summary of performance test with Neumann condition. }
	\begin{tabular}{lccc}
		$\theta$  & \multicolumn{1}{l}{$\max|y''-f|$} & \multicolumn{1}{l}{ $\max(|y_{opt}-y_{b}|)$ } & \multicolumn{1}{l}{ $\max(|y_{rk4}-y_{b}|)$ } \\ \hline\hline
		$\pi/2$ &  1.1E-07 & 8.8E-10 & 2.4E-07 \\
			$3\pi/2$ & 1.1E-06 & 1.8E-08 & 1.1E-04 \\ \hline
	\end{tabular}%
	\label{table:Neumann}
\end{table}%
\subsection{The performance with Dirichlet condition}\label{subsec:Dirichlet}
\begin{figure}[htbp]\label{fig:Dirichlet}
	\centering
	\includegraphics[width=6.4cm]{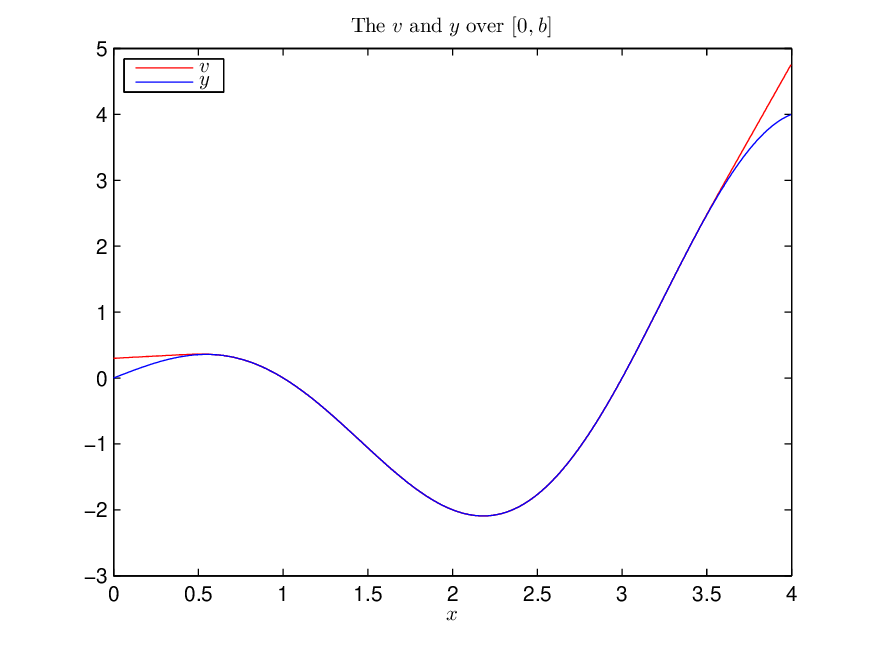}
	\includegraphics[width=6.4cm]{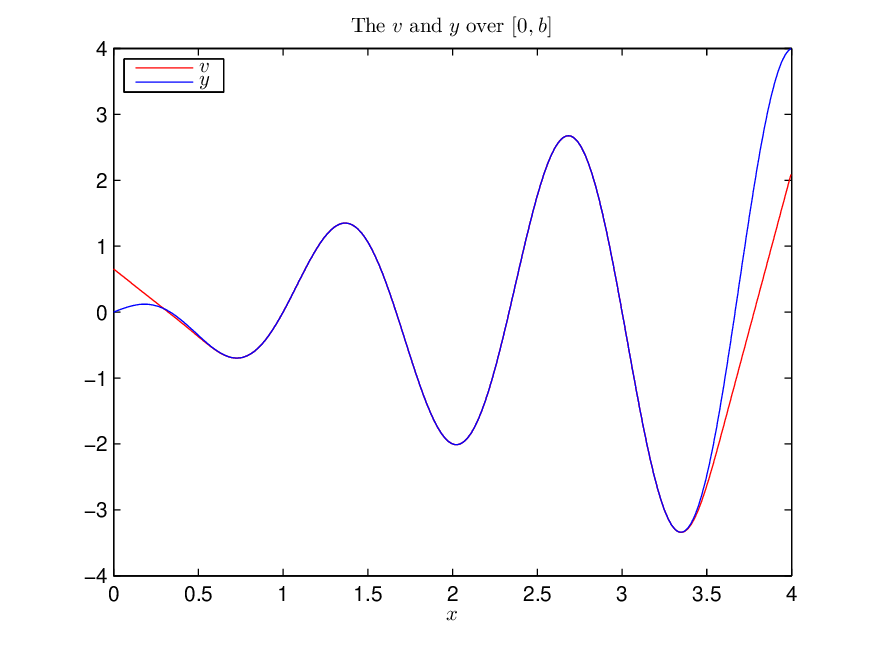}
	\includegraphics[width=6.4cm]{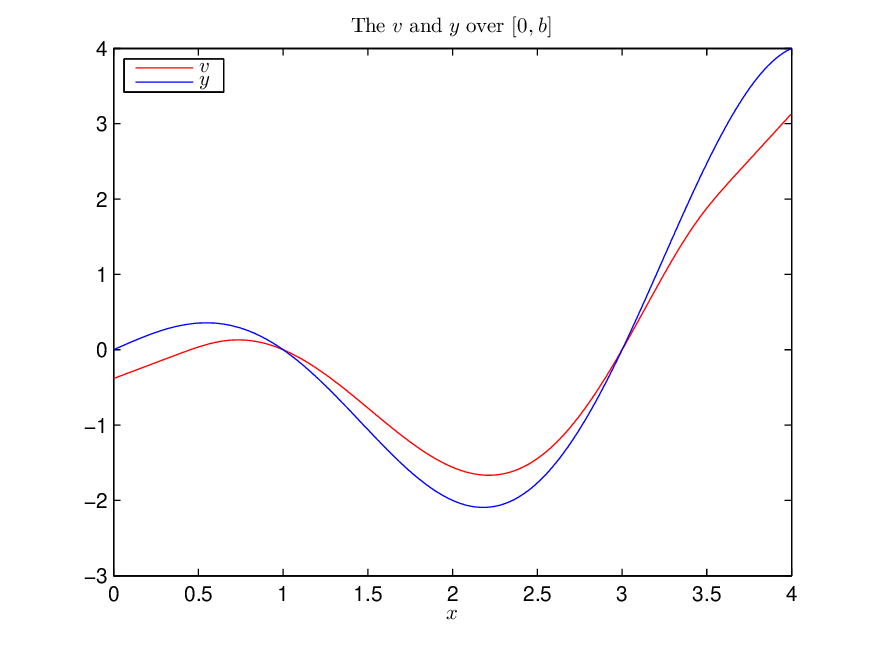}
	\includegraphics[width=6.4cm]{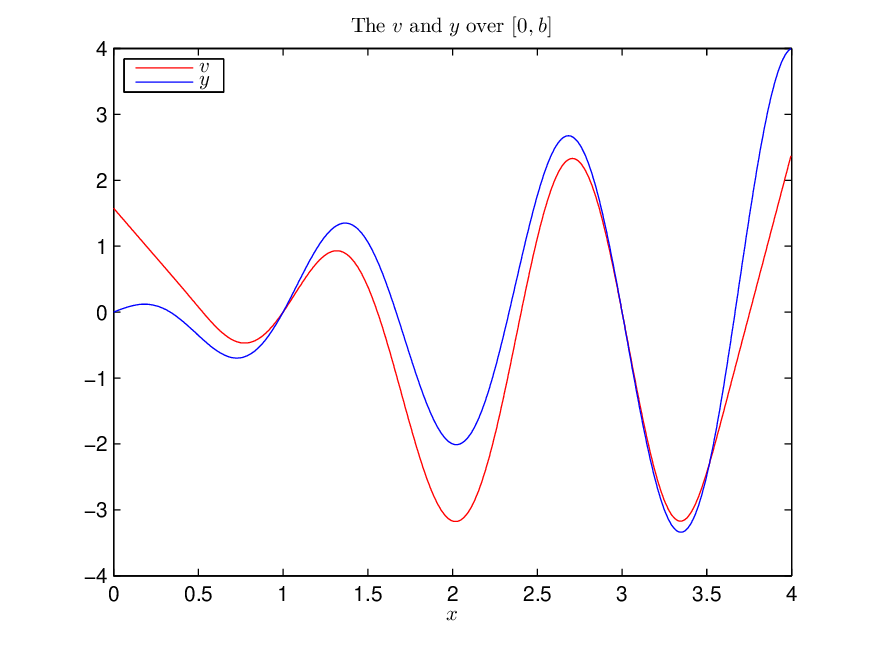}
	\caption{The plots of identified solutions (in red) as well as base function $y_b$ over $[0,b]$.  Top row shows two solutions (in red) that recover the associated base $y$ (in blue) for $\theta=\pi/2$ (top left) and  $\theta=3\pi/2$ (top right) respectively. Bottom row shows other two solutions (in red) and compares them to base $y$  for $\theta=\pi/2$ (top left) and  $\theta=3\pi/2$ (top right) respectively.  They are associated to scenarios $2$ and $1$ (top left, bottom left) for $\theta=\pi/2$ and scenarios $5$ and $1$ (top right, bottom right) for $\theta=3\pi/2$.}
\end{figure}

Table \ref{table:Dirichlet_theta2} and \ref{table:Dirichlet_theta5} summarizes the test results for $\theta=\pi/2$ and $3\pi/2$ respectively.
The optimization identifies two solutions in both $\theta=\pi/2$ and $3\pi/2$ as shown in Figure \ref{fig:Dirichlet},  one of two solutions recover $y_b$ over $[s,e]$ as expected.  In both cases,  $y_{opt}$ approaches to $y_b$ much more closer to $y_{rk4}$ for those scenarios with status $\to y_b$.  There are three divergent scenarios in total, all of them are in group $5$ where initial guess of $y(s),y'(s))$ used in shooting method is significantly far away from the true value, suggesting certain robustness of TIBO. 
\begin{table}[htbp]
	\centering
	\small
	\caption{The summary of performance test with Dirichlet condition with $\theta=\pi/2$. Severn scenarios (1,5,11-12,14,22,24) approach to $y_b$, eighteen scenarios (2-4,6-10,13,15-21,25) go to the other solution $y_s$, one scenario (23) diverges.}
	\begin{tabular}{lcccc}
		status & number    & \multicolumn{1}{l} {$\max|y''-f|$} & \multicolumn{1}{l} {$\max|y_{opt}-y_{b}|$}  & \multicolumn{1}{l} {$\max|y_{rk4}-y_{b}|$} \\ \hline\hline
		$\to y_b$   & 17  & 1.0E-07 & 4.1E-10 & 8.1E-07  \\
		$\to y_s$  & 7  & 1.1E-07 & 4.4E-01 & 2.2E-01 \\
		$diverge$ & 1 & diverge  & diverge  & diverge  \\ \hline
	\end{tabular}%
	\label{table:Dirichlet_theta2}%
\end{table}%
\begin{table}[htbp]
	\centering
	\small
	\caption{The summary of performance test with Dirichlet condition with $\theta=3\pi/2$.  Three scenarios (2,4,5)  approach to $y_b$,  nineteen scenarios (1,3,6-22) go to the other solution $y_s$, three scenarios (23-25) diverge.}
	\begin{tabular}{lcccc}
		status & number  & \multicolumn{1}{l} {$\max|y''-f|$} & \multicolumn{1}{l} {$\max|y_{opt}-y_{b}|$}  & \multicolumn{1}{l} {$\max|y_{rk4}-y_{b}|$}  \\ \hline\hline
		$\to y_b$   & 3     & 1.1E-06 & 2.6E-10 & 4.4E-05  \\
		$\to y_s$  & 19     & 1.1E-06 & 1.2E+00 & 1.2E+00  \\
	$diverge$ & 3& 	  diverge & diverge & diverge \\ \hline
	\end{tabular}%
	\label{table:Dirichlet_theta5}%
\end{table}%

 If there are more than one solutions for an ODE system, additional info can be helpful for TIBO to identify a target solution. As an example,  to force $y_{opt}\to y_b$, we pose a condition on $y'(s)$'s range such that $|y'(s)-f'(s)| \le 0.1 \cdot |f'(s)|$ for both $\theta=\pi/2$ and $\theta=3\pi/2$.  In practice, the thresholds $ f'(s)\pm 0.1 \cdot |f'(s)|$ should be based on some expectations on target solution. Table \ref{tab:Dirichlet_condition_impact_theta2} and \ref{tab:Dirichlet_condition_impact_theta5} provide updated results for $\theta=\pi/2$ and $\theta=3\pi/2$ respectively. We can see that $y_{opt}$ converges only to $y_b$ although more divergent cases are observed compared to unconditional optimization, mainly in case $\theta=3\pi/2$.
\begin{table}[htbp]
	\centering
	\small
	\caption{The impact of optimization conditioning on $|u(s)-y'(s)| \le 0.1\cdot |y'(s)|$ with $\theta=\pi/2$.
	Seventeen scenarios (2-22,24-25) approach to $y_b$, no scenarios go to the other solution $y_s$, two scenario (1,23) diverges.}
	\begin{tabular}{lccc}
		status & number  & \multicolumn{1}{l}{$\max(|y''-f|)$} & \multicolumn{1}{l}{$\max(|y_{opt}-y_{b}|)$}  \\ \hline\hline
		$\to y_b$  & 23     & 1.0E-07 & 4.1E-10  \\
		$\to y_s$   & 0     & NA & NA   \\  
		$diverge$ & 2     & diverge & diverge  \\ \hline
	\end{tabular}%
	\label{tab:Dirichlet_condition_impact_theta2}%
\end{table}%
\begin{table}[htbp]
	\centering
	\caption{The impact of optimization conditioning on $|u(s)-y'(s)| \le 0.1\cdot|y'(s)|$ with $\theta=3\pi/2$.
	Thirteen scenarios (1-6, 8-10,16,18,19,22) approach to $y_b$, no scenarios go to the other solution $y_s$, twelve scenario (7,11-15,17,20-21,23-25) diverges.
	}
	\begin{tabular}{lccc}
			status & number  & \multicolumn{1}{l}{$\max|y''-f|$} & \multicolumn{1}{l}{ $\max(|y_{opt}-y_{b}|)$} \\ \hline\hline
		$\to y_b$  & 13	     & 1.1E-06 & 2.9E-10 \\
			$\to y_s$   & 0     & NA & NA   \\  
			$diverge$ & 12     & diverge & diverge  \\ \hline
	\end{tabular}%
	\label{tab:Dirichlet_condition_impact_theta5}%
\end{table}%
\subsection{The performance with $Mix$}\label{subsec:Mix_1}
Table \ref{tab:Mix_1_summary_theta2} and \ref{tab:Mix_1_summary_theta5} report similar results under $Mix$ type for $\theta=\frac{\pi}2$ and $\theta=\frac{3\pi}2$ respectively.  

\begin{figure}[htbp]\label{fig:Mix_1}
	\centering
	\includegraphics[width=6.4cm]{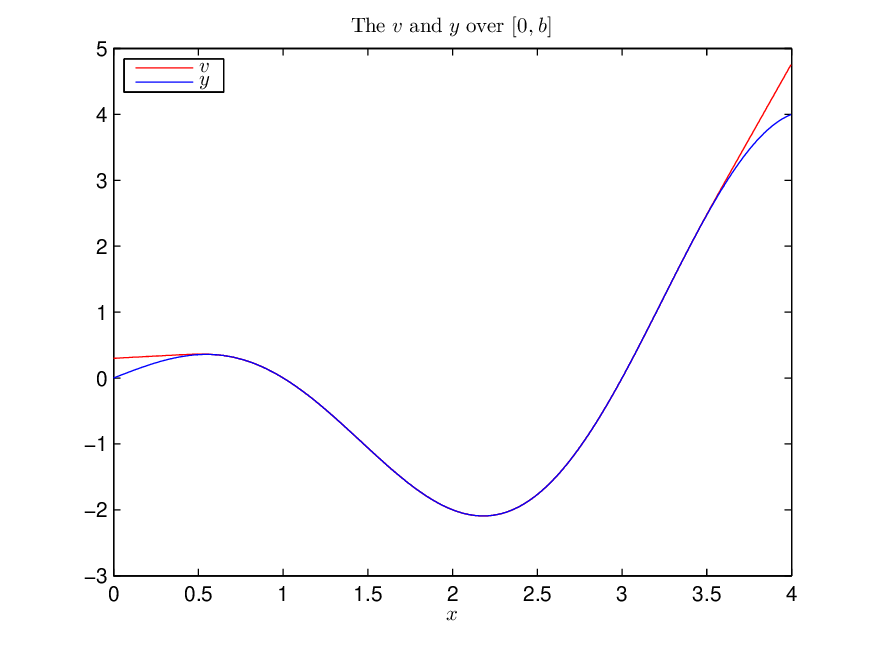}
	\includegraphics[width=6.4cm]{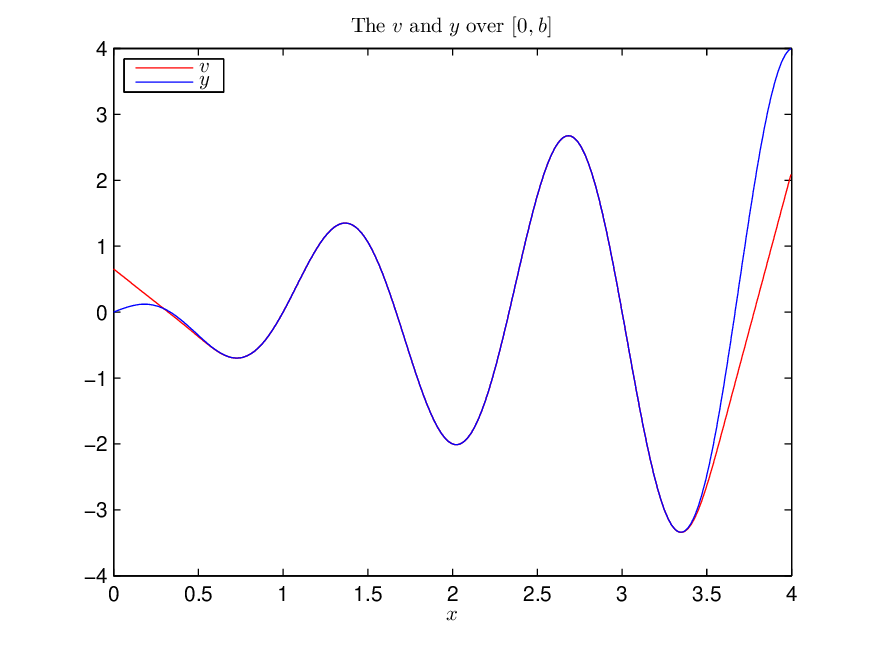}
	\includegraphics[width=6.4cm]{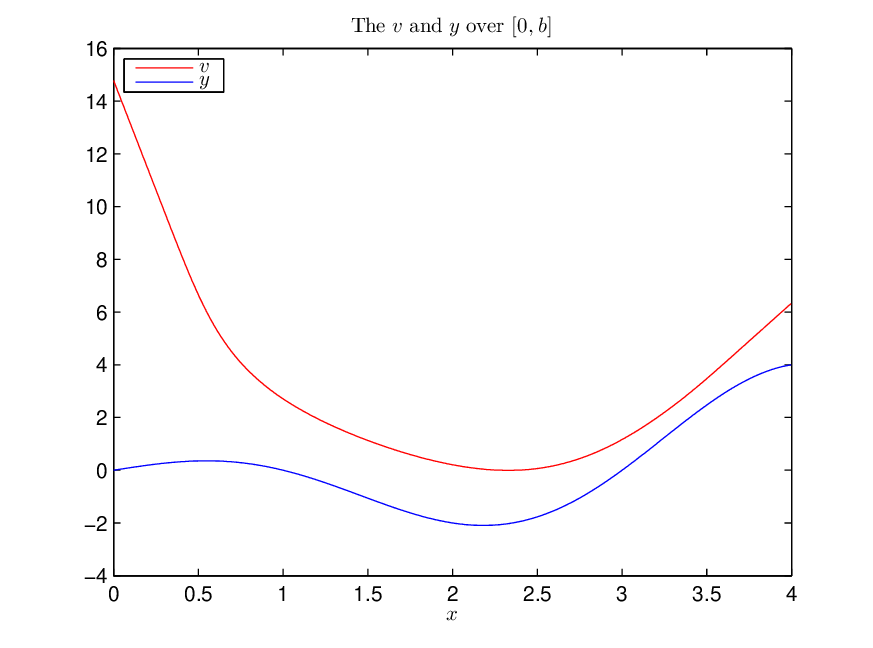}
	\includegraphics[width=6.4cm]{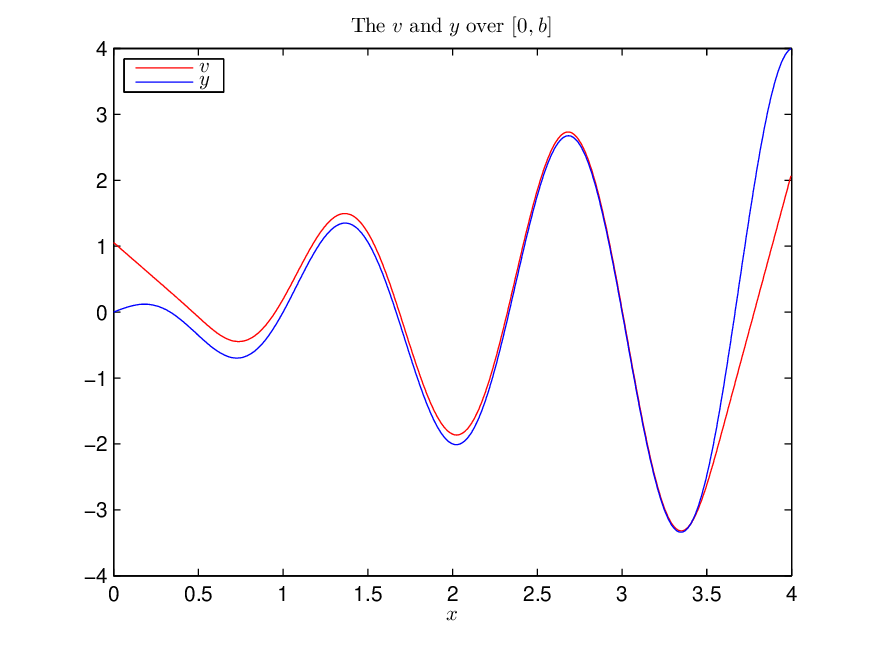}
	\caption{The plots of identified solutions (in red) as well as base function $y_b$ over $[0,b]$.  Top row shows two solutions (in red) that recover the associated base $y$ (in blue) for $\theta=\pi/2$ (top left) and  $\theta=3\pi/2$ (top right) respectively. Bottom row shows other two solutions (in red) and compares them to base $y$  for $\theta=\pi/2$ (top left) and  $\theta=3\pi/2$ (top right) respectively.  They are associated to scenarios $2$ and $1$ (top left, bottom left) for $\theta=\pi/2$ and scenarios $5$ and $1$ (top right, bottom right) for $\theta=3\pi/2$.}
\end{figure}

Same as under $Dirichlet$ type, TIBO identifies two solutions with both $\theta=\pi/2$ and $3\pi/2$ as shown in Figure \ref{fig:Mix_1}.  Compared to $Dirichlet$ type,  more scenarios diverges, which is expected since shooting method need search in $2$-dim space to identity $y(s), y'(s)$.  Nevertheless,  $y_{opt}$ converges in an acceptable rate, especially with $\theta=\pi/2$.  For all convergent cases, $y_{opt}$ keeps roughly same accuracy as under $Dirichlet$ type. On the other hand, $y_{rk4}$ performance turns to be not dependable with significant errors such that $\max|y_{rk4}-y_{b}|>10^{-2}$ for all scenarios with both $\theta=\pi/2$ and $3\pi/2$.
Figure \ref{fig:Mix_1} plots the identified two solutions $v$ and compare them to the base function $y$ with $\theta=3\pi/2$.
\begin{table}[htbp]
	\centering
	\caption{The summary of performance test under $Mix$ type with $\theta=\pi/2$.  
		Twelve scenarios (1-6,8,10,12,15,21,24) approach to $y_b$, eleven scenarios (7,9,11,13-14,16-20,23) go to the other solution $y_s$, two scenarios (22,25) diverges.  Note that $\max|y_{rk4}-y_{b}|$ values are not presented and they   are larger than $10^{-2}$}
	\begin{tabular}{lccc}
		status & number  & \multicolumn{1}{l}{$\max|y''-f|$} & \multicolumn{1}{l}{ $\max|y_{opt}-y_{b}|$ }  \\ \hline\hline
   $\to y_b$  & 12	      & 1.0E-07 & 1.3E-09   \\
  $\to y_s$  & 11	   & 1.6E-07 & 2.7E+00   \\ 
 $diverge$ & 2  & diverge & diverge   \\ \hline
	\end{tabular}%
	\label{tab:Mix_1_summary_theta2}%
\end{table}%
\begin{table}[htbp]
	\centering
	\caption{The summary of performance test under $Mix$ type with $\theta=3\pi/2$. Eight scenarios (1,3,5,6,8,9,11,14) approach to $y_b$, three scenarios (2,4,7 ) go to the other solution $y_s$,  fourteen scenarios (10,12-13,15-25) diverges. Note that $\max|y_{rk4}-y_{b}|$ values are not presented and they are all larger than $10^{-2}$.}
	\begin{tabular}{lccc}
	status & number  & \multicolumn{1}{l}{$\max|y''-f|$} & \multicolumn{1}{l}{ $\max|y_{opt}-y_{b}|$ }  \\ \hline\hline
		 $\to y_b$  & 8	     & 1.1E-06 & 6.8E-08  \\
		$\to y_s$  & 3     & 1.1E-06 & 2.0E-01  \\ 
		$diverge$ & 14   & diverge & diverge   \\ \hline
	 \hline
	\end{tabular}%
	\label{tab:Mix_1_summary_theta5}%
\end{table}%

We try another way to apply extra info of a target solution to copy with the situation where there are more than one solutions.  Note that, for the case $\theta=\pi/2$ in Figure \ref{fig:Mix_1},  low boundary of second solution $y_s$ is slightly less than $0$.  TIBO is enhanced by using conditional optimization on $y\ge -0.01$ and  the results are shown in Table \ref{tab:Mix_1_summary_theta5_codition_lower_boundary}. Except for four divergent scenarios,  all others converge to the target solution $y_s$ with decent accuracy.  
\begin{table}[htbp]
	\centering
	\small
	\caption{Impact on conditioning on lower boundary of solution with $\theta=\pi/2$.  ${y_s}$ denote the second solution identified in Table \ref{tab:Mix_1_summary_theta2}. None approaches to $y_b$, twenty one scenarios (1-2,4-14, 16-21,23-24) go to the other solution $y_s$,  four scenarios (3,15,22,25) diverges. }
	\begin{tabular}{lccc}
		status & number &  \multicolumn{1}{l}{$\max|y''-f|$}  & \multicolumn{1}{l}{ $\max|y_{opt}-{y_s}|$ } \\\hline\hline
		 $\to y_b$  & 0     & NA & NA    \\
		$\to y_s$  & 21     & 1.2E-07 &  3.1E-10 \\
		$diverge$    & 4 & diverge & diverge  \\ \hline
	\end{tabular}%
	\label{tab:Mix_1_summary_theta5_codition_lower_boundary}%
\end{table}%

\section{Summary}\label{sec:summary}
In this paper, we propose a trigonometric-interpolation based optimization (TIBO) algorithm to approximate solutions of second order nonlinear ODEs with general two-point linear boundary conditions.  The algorithm has a few advantages. 
It solves second order nonlinear ODE system with general linear boundary conditions and the idea can be extended to high order nonlinear ODE \cite{zou_tri_V}.  The algorithm is expected to converge quickly with high order of accuracy if underlying solution is sufficient smooth.  In the case that solution is not unique, TIBO provides flexibility to identify a desired solution by applying conditional optimization on certain requirements of the target solution as shown in Subsection \ref{subsec:Dirichlet} and \ref{subsec:Mix_1}.  Model performance has been tested under three types of boundary conditions.  For standard Neumann and Dirichlet condition, TIBO converges with expected accuracy and outperform significantly a benchmark method based on combination of RK4 and shooting method.  For the mixed boundary conditions, the algorithm becomes more sensitive to the initial guess of $y(s), y'(s)$ used in the optimization process and more divergent scenarios are observed. Nevertheless,  overall convergent rate is still acceptable, especially when target solutions are not highly volatile. The test results demonstrate that TIBO converges with certain robustness and tolerate certain errors of initial values used by optimization.
%\section*{Acknowledgments}
%This paper is

\newpage
\bibliographystyle{siamplain}
\bibliography{references_tibo}
\end{document}